\documentclass[times,10pt,twocolumn]{article}
\usepackage{times}
\usepackage{amsmath,amssymb,latexsym,verbatim}
\newcommand\nc{\newcommand} \nc\rnc{\renewcommand}
\nc\nn\newenvironment       \nc\nt{\newtheorem}

\nt{thm}{Theorem}[section] 
\nt{cor} [thm]{Corollary}
\nt{la}  [thm]{Lemma}
\nt{prop}[thm]{Proposition} 
\nt{prov}[thm]{Proviso} 

\nt{dftemp}  [thm]{Definition} 
\nt{rmktemp} [thm]{Remark} 
\nn{df}{\begin{dftemp}\normalfont}{$\square$\end{dftemp}} 
\nn{dfSansSquare}{\begin{dftemp}\normalfont}{\end{dftemp}} 
\nn{rmk}{\begin{rmktemp}\normalfont}{$\square$\end{rmktemp}}

\nn{pf}{\smallskip\noindent\emph{Proof}\quad}{\hfill$\square$\par\smallskip}
\nn{pfSansSquare}{\smallskip\noindent\emph{Proof}\quad}{\smallskip}

\nc\Atom{\mbox{Atom}}
\nc\B{\ger{B}}
\nc\bx{\mbox{Box}}
\nc\EF{\mbox{EF}}
\nc\ep{\varepsilon}
\nc\Orphan{\mbox{Orphan}}
\nc\Parent{\mbox{Parent}}
\nc{\U}{\Upsilon}
\nc{\arrow}{\longrightarrow}
\nc{\ger}[1]{\ensuremath{\mathfrak {#1}}}
\rnc{\phi}{\varphi}

\begin{document}
\title{Spectra of Monadic Second-Order Formulas\\ with One Unary Function}
\author{Yuri Gurevich\\ Microsoft Research\\ 
One Microsoft Way\\ Redmond, WA \ 98052
\and 
Saharon Shelah\thanks{Publication 536, 
partially supported by the US--Israel Binational Science Foundation.}\\  
Mathematics, Hebrew University\\ 
Givat Ram, 91904 Jerusalem\\ 
and Math. Dept, Rutgers University\\ 
New Brunswick, NJ \ 08903}
\date{}
\maketitle
\thispagestyle{empty}

\begin{abstract}
We establish the eventual periodicity of the spectrum of any monadic
second-order formula where\\ (i)~all relation symbols, except equality,
are unary, and\\ (ii)~there is only one function symbol and that symbol
is unary.
\end{abstract}

\section{Introduction}\label{sec:intro}

Durand, Fagin and Loescher established the eventual periodicity of the
spectrum of any monadic first-order formula with one unary function symbol
{\cite{DFL97}}.  (They also discuss the cases of two unary function
symbols and one binary function symbol; in either of those cases the
spectra are more complicated.)  In this paper, we are interested in
monadic second-order logic (MSO).

Let us recall the relevant definitions.  The {\em spectrum} of a
first-order or second-order formula is the set of the cardinalities of its
finite models.  A set $S$ of natural numbers is {\em eventually periodic}
if there exist natural numbers $\theta$ and $p > 0$ such that for every
$n\ge\theta$, if $S$ contains $n$ then it contains $n + p$.

{\em Monadic first-order formulas with one unary function symbol} are
first-order formulas subject to the following two restrictions:
\begin{itemize}
\item
All predicate symbols, with the exception of equality, are unary.  
\item
There is only one function symbol.  It is unary.
\end{itemize}
{\em MSO formulas with one unary function symbol} are defined similarly.

Upon learning the result of Durand, Fagin and Loescher, we noticed that
the monadic second-order composition method gives the following
generalization.

\begin{thm}[Main Theorem] 
The spectrum of any monadic second-order formula with one unary function
symbol is eventually periodic.
\end{thm}

Durand, Fagin and Loescher kindly mentioned our generalization at the end
of their paper.  For no reason, a proof sketch of the generalization was
lying idle all these years.  In this paper we prove the Main Theorem.
Additional results will appear in {\cite{S03}}.  We give also a direct
proof of the decidability of the finite satisfiability of MSO formulas
with one unary function symbol.  The fact of decidability is known
{\cite{R69,BGG97}}.
 
The paper is organized as follows.  Section~\ref{sec:cm} presents the
basics of the composition method.
Sections~\ref{sec:ep}--\ref{sec:op} prepare the ground for the
proof of the Main Theorem; the proof itself is given in
Section~\ref{sec:proof}.  Specifically, Sections~\ref{sec:ep} gives simple
facts on eventual periodicity, Section~\ref{sec:struc} introduces
structures of relevance to the proof of the Main Theorem,
Section~\ref{sec:theories} reduces the Main Theorem to a similar theorem
where the role of formulas is played by special finite fragments of the
theories of function graphs, and Section~\ref{sec:op} introduces
relevant operations on structures of relevance.  The direct decidability
result is proven in Section~\ref{sec:fs}.

\section{The Composition Method}\label{sec:cm}

Under certain circumstances, a composition of structures gives rise to a
composition of their appropriately defined types.  That observation lies in
the heart of the composition method.  We explain a simple version of the
method that goes a long way and that is sufficient for our purposes in
this paper.  At the end, we mention some further developments.  Lower your
expectation of the feasibility of algorithms. 

Let $L$ be a purely relational language of finite vocabulary, $X,Y$ be
models for $L$, and $x_1,x_2,\ldots$ (resp.\ $y_1,y_2,\ldots$\/) be
elements of $X$ (resp.\ elements of $Y$\/).  We write
$\phi(v_1,\ldots,v_j)$ to indicate that $\phi$ is a formula with free
variables among $v_1,\ldots,v_j$.

\begin{df}\label{df:dtheory}
The {\em $0$-theory\/} $th^0(X,x_1,\ldots,x_j)$ is the set of atomic
formulas $\phi(v_1,\ldots,v_j)$ such that $X\models\phi(x_1,\ldots,x_j)$,
and the {\em $(d+1)$-theory}\mbox{\ } $th^{d+1}(X,x_1,\ldots,x_j)$ is the
set of $d$-theories\/ $th^d(X,x_1,\ldots,x_j,x_{j+1})$ where $x_{j+1}$
ranges over the elements of $X$.
\end{df}

Let $\Atom(j)$ be the set of atomic formulas with variables among
$v_1,\ldots, v_j$, $\bx^0(j)$ be the powerset of $\Atom(j)$, and
$\bx^{d+1}(j)$ be the powerset of $\bx^d(j+1)$.  Every $\bx^d(j)$ is
hereditarily finite.

\begin{la}\label{la:finite}
Every $d$-theory $th^d(X,x_1,\ldots,x_j)$ belongs to $\bx^d(j)$.
\end{la}

\begin{pf}
Induction on $d$.
\end{pf}

It follows that every $th^d(X,x_1,\ldots,x_j)$ is finite.  A bound on the
cardinality of $th^d(X,x_1,\ldots,x_j)$ is computable from $d$ and $j$.

The $d$-theory $th^d(X,x_1,\ldots,x_j)$ is closely associated with the set
of formulas $\phi(v_1,\ldots,v_j)$ of depth $d$ such that
$X\models\phi(x_1,\ldots,x_j)$. 

\begin{la}\label{la:epsilon}
For every $t\in\bx^d(j)$, there is a formula $\ep_t(v_1,\ldots,v_j)$ of
quantifier depth $d$ such that the following holds for all $X$ and all
$x_1,\ldots,x_j$.
\begin{enumerate}
\item
$th^d(X,x_1,\ldots,x_j) = t$ if and only if $X\models\ep_t(x_1,\ldots,x_j)$.
\item 
If $th^d(X,x_1,\ldots,x_j) = t$ then, for every formula
$\phi =\phi(v_1,\ldots,v_j)$ of quantifier depth $\le d$, we have
\begin{itemize}
\item 
if $X\models\phi(x_1,\ldots,x_j)$ then $\ep_t$ implies $\phi$, and
\item
if $X\models\neg\phi(x_1,\ldots,x_j)$ then $\ep_t$ implies $\neg\phi$.
\end{itemize}
\end{enumerate}
\end{la}

\begin{pf}
Both claims are proved by induction on $d$, and in both cases the cases
$d=0$ is obvious.  Let $d = c+1 > 0$.  

\smallskip
{\em 1.}  The desired
$\ep_t(v_1,\ldots,v_j)$ is the conjunction of formulas
$\exists v_{j+1}\ep_s(v_1,\ldots,v_{j+1})$ where $s\in t$ and formulas
$\neg\exists v_{j+1}\ep_s(v_1,\ldots,v_{j+1})$ where $s\in\bx^c(j+1) - t$.  

\smallskip
{\em 2.}  Assume that $th^d(X,x_1,\ldots,x_j) = t$.  Every formula $\phi
=\phi(v_1,\ldots,v_j)$ of quantifier depth $\le d$ is equivalent to a
Boolean combination of formulas $\exists v_{j+1}\phi_i(v_1,\dots,v_{j+1})$
where the quantifier depth of every $\phi_i$ is $\le c$.  Therefore it
suffices to prove the second claim for the case when $\phi$ has the form
$\exists v_{j+1}\psi(v_1,\dots,v_{j+1})$.

First suppose that $X\models\phi(x_1,\ldots,x_j)$.  Then there exists
$x_{j+1}$ such that $X\models\psi(x_1,\ldots,x_{j+1}$.  Let $s=
th^c(X,x_1,\ldots,x_{j+1})$.  Clearly, $s\in t$.  By the induction
hypothesis, $\ep_s$ implies $\psi$.  It follows that $\exists
v_{j+1}\ep_s(v_1,\ldots,v_{j+1})$ implies $\phi$, and therefore $\ep_t$
implies $\phi$.

Second suppose that $X\models\neg\phi(x_1,\ldots,x_j)$.  By contradiction
suppose that $\ep_t$ is consistent with $\phi$.  Then there exist $Y$ and
$y_1,\ldots,y_{j+1}$ such that
$Y\models\ep_t(y_1,\ldots,y_j)\wedge\psi(y_1,\ldots,y_{j+1})$.  By the
first claim, $th^d(Y,y_1,\ldots,y_j) = t$, and so the $c$-theory $s =
th^c(Y,y_1,\ldots,y_{j+1})$ belongs to $t$.  By the induction hypothesis,
$\ep_s$ implies $\psi$.  Hence $\exists v_{j+1}\ep_s$ implies $\phi$.
Since $s\in t = th^d(X,x_1,\dots,x_j)$, we have
$X\models\ep_s(x_1,\ldots,x_j,x_{j+1})$ for some $x_{j+1}$.  Hence
$X\models\phi(x_1,\ldots,x_j)$ which is impossible.
\end{pf}

\begin{cor}\label{cor:decides}
There is an algorithm that, given a $d$-theory $th^d(X)$ and a sentence
$\phi$ of quantifier depth $d$ decides whether $X\models\phi$.
\end{cor}

\begin{cor}\label{cor:u}
For every sentence $\phi$ of quantifier depth $d$, there are $d$-theories
$t_1,\ldots,t_k$ such that any $X\models\phi$ if and only if
$th^d(X)\in\{t_1,\ldots,t_k\}$.
\end{cor}

There is a close connection between finite theories and
Ehrenfeucht-Fra\"{\i}ss\'e games.

\begin{la}\label{la:games}
The following are equivalent:
\begin{enumerate}
\item
$th^d(X,x_1,\ldots,x_j) = th^d(Y,y_1,\ldots,y_j)$,
\item
the duplicator has a winning strategy in
$EF\,^d((X,x_1,\ldots,x_j),(Y,y_1,\ldots,y_j))$.
\end{enumerate}
\end{la}

\begin{pf}
Induction on $d$.  The case $d = 0$ is obvious.  Suppose that $d = c+1 >
0$.

\smallskip
{\em 1$\rightarrow$2.}  Without loss of generality, on the first step, the
spoiler chooses some element $x_{j+1}$ in $X$.  The $c$-theory
$th^c(x_1,\ldots,x_{j+1})$ belongs to $th^d(X,x_1,\ldots,x_j)$ and
therefore to $th^d(Y,y_1,\ldots,y_j)$.  So the duplicator can find an
element $y_{j+1}\in Y$ such that $th^c(X,x_1,\ldots,x_{j+1}) =
th^c(Y,y_1,\ldots,y_{j+1})$.  In the rest of the game, the duplicator uses
the winning strategy guaranteed by the induction hypothesis.

\smallskip
{\em 2$\rightarrow$1.}  By the virtue of symmetry, it suffices to prove
that every member $th^c(X,x_1,\ldots,x_{j+1})$ of $th^d(X,x_1,\ldots,x_j)$
equals to some member $th^c(Y,y_1,\ldots,y_{j+1}))$ of
$th^d(Y,y_1,\ldots,y_j))$.  Consider a play where the spoiler starts by
choosing an element $x_{j+1}$ in $X$ and the duplicator chooses some
$y_{j+1}$ in $Y$ such that he has a winning strategy in the remainder game
$EF^d((X,x_1,\ldots,x_{j+1}),(Y,y_1,\ldots,y_{j+1}))$.  By the induction
hypothesis, $th^c(X,x_1,\ldots,x_{j+1}) = th^c(Y,y_1,\ldots,y_{j+1}))$.
\end{pf}
  
To illustrate this most basic version of the composition method, we give a
simple example.  Later in this paper we define a finite forest as a finite
acyclic directed graph, where every vertex has at most one outgoing edge,
together with a unary relation Orphan that consists of the vertices
without outgoing edges, the {\em orphans}\/.  (Actually we will use finite
forests that may be endowed with additional unary relations, but for now
let's forgo that expansion.)

Every finite forest is obtained from a singleton forest by means of the
disjoint union operation $X+Y$ and the unary operation $X'$ that attaches
a root to the given forest $X$ and thus turns it into a tree.  

\begin{la}\label{la:focm}\mbox{}
\begin{enumerate}
\item
The $d$-theory of $th^d(X+Y)$ is uniquely determined by the
$d$-theories $th^d(X)$ and $th^d(Y)$.
\item
The $d$-theory of $th^d(X')$ is uniquely determined by the
$d$-theory $th^d(X)$.
\end{enumerate}
\end{la}

\begin{pf}
Use Ehrenfeucht-Fra\"{\i}ss\'e games.
\end{pf}

This lemma can be called the first-order composition lemma for finite
forests.  It allows us to define the operation $th^d(X) + th^d(Y) =
th^d(X+Y)$ and $(th^d(X))' = th^d(X')$ on $d$-theories.  Both operations
on $d$-theories are computable.

\begin{prop}\label{prop:decidable}\mbox{}
The first-order theory of finite forests is decidable.
\end{prop}

\begin{pf}
By Corollary~\ref{cor:decides}, it suffices to show that the set $S_d$ of
the $d$-theories of finite forests is computable.  Every member of $S_d$
is obtained from the $d$-theory of the singleton forest by means of the
two operations on $d$-theories.  Since $S_d$ is finite, this allows us to
compute $S_d$.
\end{pf}

Now we turn to MSO logic.  It is convenient to represent second-order
structures as special first-order structures.

\begin{dfSansSquare}\label{df:boolassoc}
The {\em Boolean associate} $\B(X)$ of a structure $X$ is the atomic
Boolean algebra whose atoms are the element of $X$ together with
\begin{itemize}
\item the relations of $X$ as relations of $\B(X)$, and
\item an additional unary relation Atomic that consists of the atoms.
\quad\hfill$\square$
\end{itemize}
\end{dfSansSquare}

\begin{df}\label{df:msodtheory}
The MSO $d$-theory $Th^d(X)$ is the first-order $d$-theory $th^d(\B(X))$. 
\end{df}

\begin{rmk}
The definition of first-order $d$-theories above had a prerequisite: the
first-order language is purely relational.  The point of the prerequisite
is to ensure that every $\bx^0(j)$ is finite.  In the presence of function
symbols, for sufficiently large $j$, we have infinitely many terms with
variables among $v_1,\ldots,v_j$ and therefore $\bx^d(j)$ is infinite.
The prerequisite can be waved for Boolean associates because, in the case
of Boolean algebras, we have only finitely many inequivalent terms.  But
then the definition of $bx^d(j)$ should be modified so that only terms in
an appropriate normal form are used.
\end{rmk}

Our example generalizes to MSO.

\begin{la}\label{la:msocl}\mbox{}
\begin{enumerate}
\item
$Th^d(X+Y)$ is uniquely determined by $Th^d(X)$ and $Th^d(Y)$.
\item
$Th^d(X')$ is uniquely determined by $Th^d(X)$.
\end{enumerate}
\end{la}

This lemma can be called the MSO composition lemma for finite forests.  It
allows us to define the operation $Th^d(X) + Th^d(Y) = Th^d(X+Y)$ and
$(Th^d(X))' = Th^d(X')$ on $d$-theories.  Both operations on MSO
$d$-theories are computable.

Proposition~\ref{prop:decidable} and its proof generalize to the MSO case
but we give an alternative proof that bounds the size of the minimal
forest with a given MSO $d$-theory.

\begin{prop}\label{prop:bounded}\mbox{}
There is a computable function $F(d)$ such that, for every $d$ and every
finite forest $X$, there is a forest $Y$ of cardinality $\le F(d)$ with
$Th^d(X) = Th^d(Y)$.
\end{prop}

\begin{pf}
Given a number $d$, build a finite sequence $S = (X_1,X_2,\ldots)$ of
finite forests as follows.  $X_1$ is the singleton forest.  Suppose that
$X_1,\ldots,X_k$ have been constructed and let $i,j$ range over
$\{1,\ldots,k\}$.  If there is a finite forest $Y$ of the form $X_i'$ or
$X_i+X_j$ such that $Th^d(Y)$ differs from any $Th^d(X_i)$, choose any
such $Y$ and set $X_{k+1} = Y$; otherwise halt.

The sequence $S$ contains every MSO the $d$-theory of every finite forest.
The length of $S$ is bounded by the cardinality of the $\bx^d(0)$ for the
first-order language of the Boolean associates of finite forests.  The
cardinality of any $X_{k+1}$ is bounded by the double of the maximal
cardinality of the forests in $X_1,\ldots,X_k$.  This gives the desired
$F(d)$.
\end{pf}

The proposition implies that the MSO theory of finite forests is
decidable.  The decision algorithm is non-elementary, but there is no
elementary algorithm for that decision problem.  This follows from
Stockmeyer's result that the decision problem for the first-order theory
of words is known to be non-elementary [S74,R02].  Words can be seen as
special trees with additional unary relations.  The order relation (on the
positions) which is a given in Stockmeyer's theorem is MSO definable from
the successor relation which is a given in our case.

\begin{rmk}\label{rmk:subtle}
One subtlety of the composition method is that it is sensitive to the
precise choice of language.  We explain this on an the example of finite
forests where, to simplify notation, we ignore the difference between
subsets of a forest $X$ and the elements of $\B(X)$.

We mentioned above that $Th^d(X')$ is determined by $Th^d(X)$.  Let $r$ be
the root of $X'$, let $R$ be the singleton forest containing $r$, let
$A_1,\ldots,A_j$ be arbitrary subsets of $X$, let $B_i = A_i\cap\{r\}$,
and let $C_i = A_i - B_i$.  More generally, $Th^d(X',A_1,\ldots,A_j)$ is
determined by $Th^d(R,B_1,\ldots,B_j)$ and $Th^d(X,C_1,\dots,C_j)$.  But
this claim of unique determination fails if we abandon the unary relation
Orphan which might have seemed to play no role until now.  Orphan is
needed already in the case $d=0$ and $j\ge2$.  Indeed, let $x$ be any
element of $X$, $A_1 = B_1 =\{r\}$ and $A_2 = C_2 =\{x\}$.  In order to
determine  whether there is an edge from $x$ to $r$, we need to know
whether $x$ is an orphan in $X$.

This little exercise shows also another subtlety of the method.  Notice
how neatly each set $A_i$ split into subsets $B_i$ and $C_i$.  It would be
less convenient to work the first-order $d$-theories
$th^d(X,x_1,\ldots,x_j)$ instead of second-order $d$-theories
$Th^d(X,A_1,\ldots,A_j)$.   Elements are not splittable.  The method fits
MSO better than first-order logic.
\end{rmk}

The origins of the composition method can be traced back to the
Feferman-Vaught article on ``the first-order properties of products of
algebraic systems'' {\cite{FV59}}.  L\"{a}uchli introduced $d$-theories in
the context of weak monadic second-order logic; he proved that the weak
MSO theory of linear order is decidable {\cite{L68}}.  (The weak MSO is
the version of MSO where second-order quantification is restricted to
finite sets.)

Shelah generalized the method to full MSO and used it in particular to
prove in a uniform way all known decidability results for the MSO theories
of various classes of linear orders {\cite{S75}}.  He introduced
$(k_1,\ldots,k_d)$-theories comprising sentences whose prenex form has $d$
blocks of quantifiers: $k_1$ quantifiers of one kind (say, existential
quantifiers), followed by $k_2$ quantifiers of the other kind, followed by
$k_3$ quantifiers of the first kind, and so on.  A $d$-theory is a
$(k_1,\ldots,k_d)$-theory where every $k_i=1$.  In the theory of linear
orders, the main composition lemma is about the addition (that is
concatenation) of linear orders.  In {\cite{S75}}, an important role was
played by two generalizations of the classical Ramsey theorem (one for
dense orders and another for uncountable well orderings) that take
advantage of the fact that, for each $d$, the $d$-theories of linear
orders together with the derived addition operation form a semigroup.

MSO logic is arguably the right paradigm for the composition method.
Composition theorems reduce the theory of a composition of structures to a
composition of their finite theories, but --- even in the case of
first-order composition theorems --- the reduction depends on the MSO
theory of the index structure.  MSO composition theorems reduce the MSO
$(k_1,\ldots,k_d)$-theory of a composition of structures to some MSO
$(\ell_1,\ldots,\ell_d)$-theory of the index structure with disjoint unary
relations $t$ where $t(i)$ means that $t$ is the $(k_1,\ldots,k_d)$-theory
of the $i^{th}$ component.

The authors used the method on numerous occasions.  Article {\cite{G78}}
lays a technical foundation for more advanced applications of the method.
Section~3 of the survey {\cite{G85}} describes the composition method
(calling it the model-theoretic decidability technique) and Section~5
mentions various applications of the method.  See also the dissertation
{\cite{Z94}} and exposition {\cite{T97}}.  A recent sophisticated use of
the composition method over finite structues is found in {\cite{S96}}.

The power of the method is under-appreciated.  Throughout the years we saw
various problems that could be solved by the composition method.  This
claim was put to the test a couple of years ago when Alexander Rabinovich
posed a conjecture to the first coauthor who insisted that the method will
confirm or refute the conjecture.  The conjecture was confirmed
{\cite{GR00}}, and Rabinovich went on to use the method {\cite{R03}}.

The present paper uses the simple form of the composition method explained
and exemplified above.

\begin{rmk}\label{rmk:automata}
One of the reviewers asked how does the composition method compare with
the automata-theoretic method.  The application domains of the two methods
intersect.  For example, the decidability of S1S, established first by
means of B\"{u}chi automata {\cite{B62}}, has also a simple
model-theoretic proof {\cite{S75,G85}}.  Rabin used automata to prove the
decidability of the MSO theory, known as S2S, of the infinite binary tree
{\cite{R69,BGG97}}.  It is not clear whether the composition theory can be
used for the purpose.  One of the consequences of Rabin's result is the
decidability of the MSO theory of rational order.  This consequence was
proved, more directly and naturally, by the composition method
{\cite{S74}}.  This alternative proof generalizes to a class of dense
linear orders (so-called short modest linear orders) {\cite{GS79}} to which
the automata method does not seem to apply.  At least on one occasion, the
two methods were used in a complimentary way {\cite{GS85}}.  The automata
method has been used much more for establishing complexity results.
Lemma~\ref{la:complexity} is one modest example of the use of the
composition method for that purpose.
\end{rmk}

The composition method is inherently model-theoretic, like that of
Ehrenfeucht-Fra\"{\i}ss\'e; it may be a useful addition to your toolbox.

\section{Eventually Periodic Sets}\label{sec:ep}

\begin{df}
A set $S$ of natural numbers is {\em eventually periodic} if there exist
natural numbers $p>0$ (a {\em period} of $S$) and $\theta$ (a {\em
$p$-threshold} for $S$) such that
$n\in S\ \mbox{implies}\ n+p\in S\ \mbox{for all }n>\theta.$\quad 
\end{df}

The definition of eventual periodicity in {\cite{DFL97}} is similar except
that ``implies'' is replaced with ``is equivalent to''.  The following
lemma shows that this alteration makes no real difference.

\begin{la}
Assume that a set $S$ is eventually periodic with a period $p$ and a
$p$-threshold $\theta$.  There exists a natural number $\theta'$ such that
$$
 n\in S\quad \mbox{is equivalent to}\quad n+p\in S
 \quad\mbox{for all }n\ge\theta'.
$$ 
\end{la}

The number $\theta'$ could be called a {\em strict $p$-threshold} for $S$.

\begin{pfSansSquare}
The lemma is trivial if $S$ is finite, so we assume that $S$ is
infinite.  For each $i = 0,\ldots,p-1$, let
$$ 
A _i =\{n\in S :\quad n\ge\theta\quad\mbox{and}\quad n = i \mod\ p  \}.
$$ 
The desired strict $p$-threshold $\theta' = \max \{\min(A_i) :
A_i\ne\emptyset\}.\qquad\square$
\end{pfSansSquare}

\begin{cor}\mbox{}\\
Every eventually periodic set $S$ has a least period, and the least
period divides any other period of $S$. 
\end{cor}

\begin{pf}
It suffices to prove that the greatest common divisor $p$ of periods $p_1,
p_2$ is a period.  Let $\theta_i$ be a strict $p_i$-threshold and $\theta
=\max(\theta_1,\theta_2)$.  We show that $\theta$ is a $p$-threshold for
$S$.

Since $p =\gcd(p_1,p_2)$, there exist integers $a_1, a_2$ such that
$a_1p_1 + a_2p_2 = p$.  Without loss of generality, $a_1 > 0$ and $a_2 <
0$.  Suppose that $n\ge\theta$ and $n\in S$.  Since $\theta_1$ is a
$p_1$-threshold, $n + a_1p_1\in S$.  Since $\theta_2$ is a strict
$p_2$-threshold, $n + p = (n + a_1p_1) - |a_2|\cdot p_2\in S$.
\end{pf}

Recall that an arithmetic progression is a set of integers of the form
$$
 \{b + jp : j = 0,1,\ldots\}.
$$
where $p>0$.  

\begin{la}[\cite{DFL97}]
A set of natural numbers is eventually periodic if and only if it is a
finite union of arithmetic progressions and singleton sets. 
\end{la}

\begin{cor}\label{cor:union}
The union of finitely many eventually periodic sets is eventually
periodic.
\end{cor}

\begin{la}
If the sets $S_1,\ldots, S_m$ are eventually periodic then so is the set
$$
 \{ n_1 + \cdots + n_m :\ \mbox{every } n_i\in S_i \}
$$
\end{la}

\begin{pfSansSquare}
Suppose that $S_i$ is eventually periodic with period $p_i$ and
$p_i$-threshold $\theta_i$, and let $p$ be the least common multiple of
$p_1,\ldots,p_m$.  Then $S$ is eventually periodic with period $p$ and a
$p$-threshold $\theta =\theta_1 +\cdots\theta_m$.  Indeed suppose that
$n\in S$ and $n\ge\theta$.  Then $n = n_1 +\cdots + n_m$ where $n_i\in S_i$.
There is an index $i$ such that $n_i\ge\theta_i$.  We have
$$
 n + p = n_1 +\cdots + n_{i-1} + (n_i + p) + n_{i+1} + \cdots + n_m
 \in S.\qquad\square
$$
\end{pfSansSquare}  

\section{Structures of Interest}\label{sec:struc}

In the rest of the paper, all structures are finite.  It is often
convenient to view unary relations as sets.

\begin{df}
A {\em partial-function graph\/} is a directed graph where every vertex
has at most one outgoing edge together with the unary relation Orphan that
consists of the vertices without outgoing edges.  The expression
``partial-function graph'' may be abbreviated to ``PF-graph''.
\end{df}

The relation Orphan is added for technical reasons; we saw already in
Remark~\ref{rmk:subtle} that the composition method is sensitive to the
language.  Orphan allows us to express the absence of orphans in a
quantifier-free way: $\Orphan =\emptyset$.

\begin{df}
A {\em function graph\/} is a PF-graph where every vertex has exactly one
outgoing edge.
\end{df}

The edge relation of a PF-graph will be denoted by $E$.  If $E(x,y)$ holds
(so that there is an edge from $x$ to $y$), then $y$ is the {\it parent}
of $x$, and $x$ is a {\it child} of $y$.  Children of the same parent are
{\it siblings}.  A vertex without a parent is an {\em orphan}.

A nonempty sequence $x_0,\ldots,x_n$ of vertices such that
$E(x_i,x_{i+1})$ holds for every $i = 0,\ldots, n-1$ is a {\em path of
length $n$ from $x_0$ to $x_n$}.  A {\em path} is {\em trivial} if its
length is zero; otherwise the path is {\em proper}.  If there is a path
from $x$ to $y$ then $x$ is a {\em descendant} of $y$ and $y$ is an {\em
ancestor} of $x$.  If there is a proper path from $x$ to $y$ then $x$ is a
{\em proper descendant} of $y$ and $y$ is a {\em proper ancestor} of
$x$.

An $n$-vertex substructure of a PF-graph is a {\em cycle} if the $n$
vertices form a proper path $x_0,\ldots, x_{n-1}, x_0$.  A vertex is {\em
cyclic} if it belongs to a {\em cycle}.  A PF-graph without a cycle is
{\em acyclic}.

A {\em colored PF-graph} $X$ is a PF-graph together with a finite
collection of unary relations, the {\em colors}.  The set of the names of
the unary relations is the {\em palette} of $X$.  It is convenient to
think of colors as sets of vertices.

\begin{prov}
In the rest of this paper, a PF-graph is always finite and colored.
\end{prov}

A PF-graph is {\em connected} if every two vertices have a common
ancestor.  A connected PF-graph has at most one orphan.

\begin{df}
A {\em tree} is a connected PF-graph with an orphan; the orphan is
the {\em root} of the tree.  
\end{df}

A connected function graph is a connected PF-graph without an orphan.  It
has a unique cycle and is formed by the cyclic vertices and their
descendants.

In any PF-graph, the relation ``$x$ and $y$ have a common ancestor'' is an
equivalence relation.  It partitions the PF-graph into connected
components.

\begin{df}
A {\em forest} is an acyclic PF-graph.
\end{df}

A PF-graph $X$ is a forest if and only if every component of $X$ is a tree
if and only if every component of $X$ has an orphan.  A PF-graph $X$ is a
function graph if and only if no component of it has an orphan.

\section{Finite Theories}\label{sec:theories}

An {\em MSO PF-graph formula} $\chi$ is an MSO formula in the
vocabulary of PF-graphs.  The {\em spectrum} of $\chi$ is the set of
the cardinalities of PF-graphs satisfying $\chi$.

\begin{la}\label{la:phi2chi}
For every MSO formula $\phi$ with one unary function
symbol, there exists an MSO PF-graph formula $\chi$ such
that (i)~every PF-graph satisfying $\chi$ is a function graph and (ii)~the
spectrum of $\phi$ equals the spectrum of $\chi$.
\end{la}

\begin{pf}\rm  
We start with a construction.  First, replace in $\phi$ the unary function
with its graph.  For example, a subformula $f(x) = y$ becomes
$E(x,y)$, and a subformula $f(f(x)) = f(y)$ may become
$$
 \exists x'\exists z(E(x,x')\land E(x',z)\land E(y,z))
$$
Second, augment the resulting PF-graph formula with a conjunct $\Orphan
=\emptyset$.  This gives the desired $\chi$.

Claim (i) is obvious.  To prove (ii), notice that every model $X$ of
$\phi$ gives rise to a function graph $Y$ satisfying $\chi$ such that $|X| =
|Y|$; in fact, the underlying set of $Y$ is that of $X$.  Let
$\phi'$ be the formula obtained from $\chi$ by replacing every atomic
formula $E(x,y)$ with equation $f(x) = y$; clearly $\phi'$ is equivalent
to $\phi$.  Every function graph $Y$ satisfying $\chi$ gives rise to a model
$X'$ for $\phi'$, and therefore for $\phi$, such that $|Y| = |X'|$; in
fact, the underlying set of $X'$ is that of $Y$.
\end{pf}

Fix an arbitrary palette $\pi$.

\begin{prov}
In the rest of the paper, PF-graphs are of palette $\pi$.  Accordingly,
MSO PF-graph formulas use only colors of palette $\pi$.
\end{prov}

According to Definition~\ref{df:boolassoc}, the {\em Boolean associate}
$\B(X)$ of a PF-graph $X$ is the finite Boolean algebra with the vertices
of $X$ as the atoms together with the edge relation of $X$ as a binary
relation of $\B(X)$, the colors of $X$ as unary relations of $\B(X)$, and
an additional unary relation Atomic that consists of the atoms.  Let $\U$
be the vocabulary of $\B(X)$.  Notice that the vocabulary $\U$ does not
depend on the choice of $X$.

\begin{la}\label{la:chi2psi}
For every MSO PF-graph formula $\chi$, there is a first-order
$\U$-formula $\psi$ such that, for every PF-graph $X$, we have
$$ 
  X\models\chi\quad\mbox{if and only if}\quad\B(X)\models\psi 
$$
\end{la}

\begin{pf}
Obvious.
\end{pf}

According to Definition~\ref{df:msodtheory}, the MSO $d$-theory $Th^d(X)$
of a PF-graph $X$ is the first-order $d$-theory $th^d(\B(X))$ of the
Boolean associate $\B(X)$ of $X$.  In the rest of the paper, we deal only
with MSO $d$-theories $Th^d(X)$ of PF-graphs $X$ and never with first-order
$d$-theories $th^d(X)$.  And so the $d$-theory of a PF-graph $X$ will mean
$Th^d(X)$.  

By Lemma~\ref{la:games}, two PF-graphs $X$ and $Y$ have the same
$d$-theory if and only if the duplicator has a winning strategy in the
$d$-step Ehrenfeucht-Fra\"{\i}ss\'e game $\EF\:^d(\B(X),\B(Y))$.  Every
$d$-theory $t = Th^d(X)$ will be called a {\em finite theory} of PF-graphs
of depth $d$.  If $X$ is a function graph then $t$ contains the formula
$\Orphan =\emptyset$ and so every PF-graph $Y$ with $Th^d(Y) = t$ is a
function graph; in this case we say that $t$ is a {\em finite theory} of
function graphs.

\begin{df}
The {\em spectrum} of a finite theory $t$ of PF-graphs of depth $d$
is the set of the cardinalities of the PF-graphs $X$ with $Th^d(X) =
t$.  The {\em spectrum} of a first-order $\U$-formula $\psi$ is the set of
the cardinalities of PF-graphs $X$ such that $\B(X)\models\psi$.
\end{df}

\begin{prop}\label{cor:finite-theory}
In order to prove the Main Theorem, it suffices to prove that the spectrum
of any finite theory of function graphs is eventually periodic.
\end{prop}

\begin{pf}
Assume that the spectrum of any finite theory of function graphs is
eventually periodic, let $\phi$ be any MSO formula with one unary function
symbol, and let $\chi$ be as in Lemma~\ref{la:phi2chi}.  Every PF-graph
satisfying $\chi$ is a function graph.  It suffices to prove that the
spectrum of $\chi$ is eventually periodic.

Let $\psi$ be as in Lemma~\ref{la:chi2psi}.  Every PF-graph satisfying
$\chi$ is a function graph.  It suffices to prove that the spectrum of
$\psi$ is eventually periodic.

Let $d$ be the quantifier depth of $\psi$.  By Corollary~\ref{cor:u},
there are $d$-theories $t_1,\ldots,t_j$ of function graphs such that any
function graph $X\models\psi$ if and only if $Th^d(X)\in\{t_1,\ldots,t_j\}$.
By Corollary~\ref{cor:union}, the spectrum of $\psi$
of eventually periodic.  
\end{pf}

\begin{rmk}
Since the definable relation Orphan is present in the vocabulary of
PF-graphs, we should be careful with the notion of a substructure of a
PF-graph $X$.  Consider the following requirement of a subset $A$ of $X$:
for every edge $E(x,y)$ of $X$, if $A$ contains $x$ then it contains $y$.
If the requirement is satisfied then the normal definition of a
substructure of $X$ generated by subset $A$ works correctly; otherwise the
normal definition does not work correctly.  We will use the notion of
substructure only in cases when the normal definition works.
\end{rmk}

\section{Composing Graphs and their Theories}\label{sec:op} 

\begin{df}\label{df:sum}
The {\em sum} $X + Y$ of PF-graphs $X, Y$ is the PF-graph that
is the disjoint union of $X$ and $Y$.
\end{df}

Every PF-graph is the sum of its connected components.  If $X = X_1
+\cdots + X_m$ and $Y = Y_1 +\cdots + Y_n$ then $X + Y = X_1 +\cdots + X_m
+ Y_1 +\cdots + Y_n$.  If $X,Y$ are function graphs then so is $X+Y$.

If $X$ and $Y$ are disjoint they will be presumed to be substructures of
$X+Y$.  If $X$ are $Y$ are not disjoint, one or both of them should be
replaced with isomorphic copies in some canonic way.  (Of course, we are
interested in structures only up to their isomorphism type, so it does not
make any difference how the isomorphic copies are chosen.  It is
convenient though to work with structures rather than with their
isomorphic types.  In particular, it is convenient that the components of
a PF-graph form a set rather than a multiset.)

\begin{la}\label{la:stable}
$Th^d(X + Y)$ is uniquely determined by $Th^d(X)$ and $Th^d(Y)$.
\end{la}

\begin{pf}
Use the Ehrenfeucht-Fra\"{\i}ss\'e games.  As far as the Duplicator is
concerned, the game $\EF\:^d(X+Y,X'+Y')$ splits into $\EF\:^d(X,X')$ and
$\EF\:^d(Y,Y')$.
\end{pf}

The lemma justifies the following operation on finite theories of the same
depth.

\begin{df}
$Th^d(X) + Th^d(Y) = Th^d(X + Y)$.
\end{df}

Obviously the operation is commutative and associative.  For technical
reasons, we introduce the notion of dotted PF-graphs.

\begin{df}
A {\em dotted PF-graph} $(X,a)$ is a PF-graph $X$ with a
distinguished element $a$, the {\em dot} of $(X,a)$.
\end{df}

The {\em Boolean associate} $\B(X,a)$ of $(X,a)$ is the Boolean associate
$\B(X)$ of $X$ expanded with a distinguished element $a$.  The MSO
$d$-theory $Th^d(X,a)$ of a dotted PF-graph $(X,a)$ is the first-order
$d$-theory $th^d(\B(X,a))$ of the Boolean associate of $(X,a)$.

\begin{dfSansSquare}\label{df:sums}
Let $(X,a)$ and $(Y,b)$ be dotted PF-graphs.
\begin{enumerate}
\item The {\em dotted sum\/} $X\dotplus (Y,b)$ is the non-dotted PF-graph
$Z$ obtained from $X+Y$ by making $b$ the parent of every orphan in $X$.
\item  The {\em dotted sum\/} $(X,a)\dotplus (Y,b)$ is the dotted PF-graph 
$(X\dotplus (Y,b),a)$.
\item The {\em circular sum\/} $(X,a)\oplus (Y,b)$ is the non-dotted
PF-graph $Z$ obtained from $X+Y$ by making $b$ the parent of every orphan
in $X$ and making $a$ the parent of every orphan in
$Y$.\quad\hfill$\square$
\end{enumerate}
\end{dfSansSquare}

\begin{la}\mbox{}
\begin{enumerate} 
\item $Th^d(X\dotplus (Y,b))$ is uniquely determined by $Th^d(X)$ and
$Th^d(Y,b)$.
\item $Th^d((X,a)\dotplus (Y,b))$ is uniquely determined by
$Th^d(X,a)$ and $Th^d(Y,b)$.
\item $Th^d((X,a)\oplus (Y,b))$ is uniquely determined by $Th^d(X,a)$ and
$Th^d(Y,b)$.
\end{enumerate}
\end{la}

\begin{pf}
Use the Ehrenfeucht-Fra\"{\i}ss\'e games.
\end{pf}

The lemma justifies the following operations on finite theories of the
same depth.

\begin{dfSansSquare}
Let $(X,a)$ and $(Y,b)$ be dotted PF-graphs.
\begin{enumerate}
\item $Th^d(X)\dotplus Th^d(Y,b) = Th^d(X \dotplus (Y,b))$
\item $Th^d(X,a)\dotplus Th^d(Y,b) = Th^d((X,a)\dotplus (Y,b))$
\item $Th^d(X,a)\oplus Th^d(Y,b) = Th^d((X,a)\oplus
(Y,b))$\quad\hfill$\square$
\end{enumerate}
\end{dfSansSquare}

The second operation is associative; we will use that fact.

\section{Proof of the  Main Theorem}\label{sec:proof}

Let $d$ be an arbitrary natural number.  For each $d$-theory $\sigma$ of
function graphs, let $p_\sigma$ be the cardinality of the smallest
function graph with $d$-theory $\sigma$.  For each $d$-theory $\tau$ of
dotted PF-graphs, let $q_\tau$ be the cardinality of the smallest
dotted PF-graph with $d$-theory $\tau$.  Let $p$ be the least common
multiple of the numbers $p_\sigma$ and $q_\tau$ for $\sigma$ and $\tau$ as
above.  We prove that the spectrum of any $d$-theory $s$ of function graphs
is eventually periodic with period $p$.  By
Proposition~\ref{cor:finite-theory}, this implies the Main Theorem.

So let $s$ be an arbitrary $d$-theory of function graphs and let $X$ be a
sufficiently large function graph with $d$-theory $s$.  We construct a
function graph $Y$ of cardinality $|Y| = |X| + p$ such that $Th^d(Y) = s$.
The meaning of sufficiently large will be clarified in the course of the
construction.  We consider several cases.

Case 1: The number of components of $X$ exceeds the number of $d$-theories
of function graphs.

Let $X_1,\ldots,X_m$ be the components of $X$ in some order, and let $Y_i$
be the sum $X_1 +\cdots + X_i$ of the first $i$ components.  We use the
fact that the sum operation on PF-graphs is associative.  By the
definition of Case 1, there exist $i<j$ such that $Th^d(Y_i) = Th^d(Y_j)$.
Let $t = Th^d(Y_i)$, $u = Th^d(X_{i+1} +\cdots + X_j)$, and $v =
Th^d(X_{j+1}+\cdots + X_m)$.  Then $t + u = t$ and $t + v = s$.
Furthermore, $t + k\cdot u = t$ for every positive $k$.  By the definition
of $p_u$, there is a function graph $Z$ such that $Th^d(Z) = u$ and $|Z| =
p_u$.  To obtain the desired $Y$, augment $X$ with $k = p/p_u$ components
isomorphic to $Z$.  Clearly, $|Y| = |X| + k\cdot p_u = |X| + p$.  The
order of the components does not matter of course but it is convenient to
imagine that the new components come after $X_j$ and before $X_{j+1}$.
Then
$$
 Th^d(Y) = t + k\cdot u + v = t + v = s
$$

Case 2: $X$ has a vertex $a$ such that the number of noncyclic children of
$a$ exceeds the number of $d$-theories of forests.  

Each noncyclic child $y$ of $a$, together with its descendants, forms a
tree; let us call it $T_y$.  Let $F$ be the forest formed by the trees
$T_y$, so that the number of the components of $F$ exceeds the number of
the $d$-theories for forests.  By the argument of Case 1, applied to
forests rather than function graphs, there exists a forest $G$ of
cardinality $|G| = |F| + p$ such that $Th^d(G) = Th^d(F)$.  

Let $Q$ be the function graph obtained from $X$ by removing the trees
$T_y$ so that $X = F\dotplus (Q,a)$.  The desired function graph $Y =
G\dotplus (Q,a)$.  

Case 3: $X$ has a cycle of cardinality $m$ that exceeds the number of
$d$-theories of dotted trees.

Arrange the vertices of the cycle into a path $a_0,\ldots,a_{m-1},a_0$.
Let $T_i$ be the tree formed by the vertex $a_i$ and its non-cyclic
descendants.  We use the fact that the dotted sum operation on dotted
PF-graphs is associative.  For each $i\in\{1,\ldots,m\}$, let $(X_i,a_0) =
(T_0,a_0)\dotplus\cdots\dotplus(T_{i-1},a_{i-1})$, so that $(X_i,a_0)$ is
a dotted tree with root $a_{i-1}$.

By the definition of Case 3, there exist $i<j$ such that $Th^d(X_i,a_0) =
Th^d(X_j,a_0)$.  Let $(U,a_i) = (T_i,a_i)\dotplus\cdots\dotplus
(T_{j-1},a_{j-1})$ and let $t = Th^d(X_i,a_0)$, $u = Th^d(U,a_i)$ so that
$t = t\dotplus u$.  Let $V$ be the PF-graph obtained from $X$ by
removing the tree $X_j$.  $(V,a_j)$ is a dotted tree with root $a_{m-1}$,
and $X = (X_j,a_0)\oplus(V,a_j)$.  Let $v = Th^d(V,a_j)$ so that $s =
t\oplus v$.

There exists a dotted tree $(Z_1,c_1)$ such that $Th^d(Z_1,c_1) = u$ and
$|Z_1| = q_u$.  Let $(Z,c)$ be the dotted sum of $p/q_u$ copies of
$(Z_1,c_1)$ so that $|Z| = p$.  The desired $Y = \big((X_j,a_0)\dotplus
(Z,c)\big)\oplus (V,a_j)$.  Clearly, $|Y| = |X| + p$.  Since $t\dotplus u
= t$ and $t\oplus v = s$, we have
$$
 Th^d(Y) = \big(t\dotplus (p/q_u)\cdot u\big) \oplus v = t\oplus v = s
$$

Case 4: $X$ has a path $a_0,\ldots,a_{m-1}$ composed of non-cyclic
elements such that $m$ exceeds the number of $d$-theories of trees.

This case is similar to Case 3.  Trees $T_i$ and dotted trees $(X_i,a_0)$
are defined as above but this time we are interested in undotted trees
$X_i$.  By the definition of Case 4, there exist $i<j$ such that
$Th^d(X_i) = Th^d(X_j)$.  Let $(U,a_i) = (T_i,a_i)\dotplus\cdots\dotplus
(T_{j-1},a_{j-1})$ and let $t = Th^d(X_i)$, $u = Th^d(U,a_i)$ so that $t =
t\dotplus u$.  Let $V$ be the function graph obtained from $X$ by removing
the tree $X_j$.  Clearly $X = X_j\dotplus(V,a_j)$.  Let $v = Th^d(V,a_j)$
so that $s = t\dotplus v$.

As in Case 3, there is a dotted tree $(Z,c)$ such that $|Z| = p$ and
$t\dotplus Th^d(Z,c) = t$.  The desired $Y = X_j\dotplus (Z,c)\dotplus
(V,a_j)$.  Clearly, $|Y| = |X| + p$.  Further,
$$
 Th^d(Y) = (t\dotplus Th^d(Z,c)) \dotplus v = t\dotplus v = s
$$
This completes Case 4.

Finally, let $K$ be the class of function graphs $X$ that do not fall into
any of the four cases.  Since $X$ does not fall into Case 1, it has only
so many components.  Since $X$ does not fall into Case 3, the cycle of any
component of $X$ is only so long.  Since $X$ does not fall into cases 2
and 4, every cyclic element of $X$ has only so many descendants.  It
follows that there is a bound $\theta$ on the cardinality of any member of
$K$.  Thus, the spectrum of $s$ is eventually periodic with period $p$ and
$p$-threshold $\theta+1$.  The main theorem is proved.

\begin{rmk}
In the original proof sketch, we used the finite Ramsey theorem in Cases 3
and 4.  Writing up the proof we realized that these cases are similar to
Case 1 where only the associativity of the sum operation was used.
\end{rmk}

\section{Finite Satisfiability}\label{sec:fs}

\begin{la}\label{la:complexity}
There exists a computable function $F(d)$ such that, for every\mbox{\ }$d$
and every function graph $X$, there is a function graph $Y$ of cardinality
$\le F(d)$ such that $Th^d(X) = Th^d(Y)$.
\end{la}

\begin{pf}
The proof is similar to the proof of Proposition~\ref{prop:bounded}.  For
the reader's convenience, we make it independent from
Section~\ref{sec:proof}.  All we have to do is to show that every function
graph can be constructed from singleton function graphs by means of the
tree operations of Section~\ref{sec:op}.  Unfortunately this is not quite
true.  We need to revise two aspects of this plan.

First, it will be easier to deal with a larger class $K$ of finite
structures: function graphs, forests, dotted forests and singleton dotted
function graphs.  

Second, we need to extend Definitions~\ref{df:sum} and~\ref{df:sums}.

\begin{dfSansSquare}
Let $(X,a)$ and $(Y,b)$ be dotted PF-graphs.
\begin{enumerate}
\item $X + (Y,b) = ((X + Y),b)$.
\item $(Y,b)\dotplus X = ((X\dotplus (Y,b)),b)$.
\end{enumerate}
\end{dfSansSquare}

For\mbox{\ }any $d$, the $d$-theories $Th^d(X + (Y,b))$ and
$Th^d((Y,b)\dotplus X)$ are uniquely determined by $Th^d(X)$ and
$Th^d(Y,b)$.  Now we can carry out our plan.

Call a $K$ structure {\em good}\/ if it is obtained from singleton $K$
structures by means of the tree operations.  It suffices to prove that all
$K$ structures are good.  By contradiction assume that there is a bad $K$
structure and let $\cal X$ be a bad $K$ structure of minimal cardinality.
Clearly $\cal X$ cannot be singleton.  To get the desired contradiction,
we show that $\cal X$ is a composition of $K$ structures of smaller
cardinality.  Clearly $X$ is connected; otherwise it is the sum of its
components.  We have three cases.

\smallskip
1.  $\cal X$ is a function graph $X$.  The cyclic elements of $X$ form a
path $a_0,\ldots,a_{m-1},a_0$.  If $m=1$, so that there is a unique cyclic
element, then $X = Y\dotplus(Z,a_0)$ where $Y$ is the forest of the
non-cyclic elements of $X$ and $Z$ is the singleton function graph
containing $a_0$.  So $m>1$.  Let $Y$ be the tree formed by $a_0$ and its
non-cyclic descendants, and let $Z$ be the remaining part of $X$.  Then $X
= (Y,a_0)\oplus (Z,a_1)$.

\smallskip
2.  $\cal X$ is a tree $X$.  Let $b$ be the root of $X$, $Y$ be the forest
of the non-root elements of $X$, and $Z$ be the singleton forest formed by
$b$.  Then $X = Y\dotplus (Z,b)$.

\smallskip
3.  $\cal X$ is a dotted tree $(X,a)$.  Let $b$ be the root of $X$, $Y$ be
the forest of the non-root elements of $X$ and $Z$ be the singleton forest
formed by $b$.  If $a = b$ then $(X,a) = (Z,a)\dotplus Y$.  Otherwise
$(X,a) = (Y,a)\dotplus (Z,b)$ where $b$ is the root of $X$.
\end{pf}

\begin{thm}
The finite satisfiability problem for MSO formulas with one unary function
symbol is decidable.
\end{thm}

\begin{pf}
Use Lemma~\ref{la:phi2chi} and the lemma above.
\end{pf}

\subsection*{Acknowledgments}

Thanks to Andreas Blass, Mike Bernadskiy, Erich Gr\"adel, Klaus Reinhardt
and the anonymous reviewers who helped us to improve the paper.

\end{document}